\documentclass[12pt]{amsart}
\usepackage{amsfonts}
\usepackage{amscd}
\usepackage{amssymb}
\usepackage{graphics}

\usepackage{amsmath}

\newcommand{\BZ}{{\mathbb{Z}}}
\newcommand{\BN}{{\mathbb{N}}}
\newcommand{\BR}{{\mathbb{R}}}

\newcommand{\gD}{\Delta}

\newcommand{\gC}{\Gamma}
\newcommand{\gc}{\gamma}
\newcommand{\gs}{\sigma}

\newcommand{\gO}{\Omega}

\newcommand{\gep}{\epsilon}

\newcommand{\ti}[1]{\tilde{#1}}

\newcommand{\SO}{\text{SO}}

\newtheorem{prop}{Proposition}[section]
\newtheorem{thm}[prop]{Theorem}
\newtheorem{lem}[prop]{Lemma}

\theoremstyle{definition}

\newtheorem{rem}[prop]{Remark}
\newtheorem{rems}[prop]{Remarks}

\newtheorem{clm}[prop]{Claim}

\begin{document}
\author{Tsachik Gelander}

\thanks{The author was partially supported by NSF grant
DMS-0404557, BSF grant 2004010, and the 'Finite Structures' Marie
Curie Host Fellowship, carried out at the Alfr\'ed R\'enyi Institute
of Mathematics in Budapest.}

\date{July 18,~2006}

\title{On deformations of $F_n$ in compact Lie groups}
\maketitle

\let\languagename\relax 

\begin{abstract}
We study some properties of the varieties of deformations of free groups in compact Lie groups. In particular we prove a conjecture of Margulis and Soifer about the density of non-virtually free points in such variety, and a conjecture of Goldman on the ergodicity of the action of $\text{Aut}(F_n)$ on such variety when $n\ge 3$. 
\end{abstract}

For $n>1$, a generic $n$--tuple of elements in a connected compact
non-abelian Lie group $G$ generates a free group. G.A. Margulis
and G.A. Soifer conjectured (c.f. \cite{Soifer}) that every such tuple
can be slightly deformed to one which generates a group which is
not virtually free. In this note we prove this conjecture,
and actually show that for $n\geq 3$ and for an arbitrary dense
subgroup $\gC$, with some restriction on the minimal size of a
generating set, the set of deformations of $F_n$ whose image coincides with
$\gC$ is dense in the variety of all deformations. 
The idea is to move a given (almost arbitrary) $n$--tuple into any open subset of $G^n$ by applying
Nielsen transformations.
Using the same idea we
also prove a conjecture of W.M. Goldman \cite{Goldman} on the
ergodicity of the action of $\text{Out}(F_n)$ on
$\text{Hom}(F_n,G)/G$ when $n\geq 3$. For $n=2$, we prove the Margulis--Soifer conjecture
by showing that any pair can be slightly deformed to one which generates an infinite group which has Serre's
property $(FA)$ and in particular is not virtually free.

\section{The case $n>2$}

We will say that a group $\gC$ is $k$--generated if it has a
generating set of size $\leq k$.

\begin{thm}\label{thm1} Let $n\geq 3$. Let $G$ be a connected compact Lie group, and let
$\gC\leq G$ be an $(n-1)$--generated dense subgroup. Then any $n$
elements $s_1,\ldots,s_n\in G$ admit an arbitrarily small
deformation $t_1,\ldots,t_n$ with $\gC=\langle
t_1,\ldots,t_n\rangle$. In other words, the set
$$
 \{ f\in\text{Hom}(F_n,G):f(F_n)=\gC\}
$$
is dense in $\text{Hom}(F_n,G)$.
\end{thm}

For example, Theorem \ref{thm1} implies that if $n$ is larger than
the minimal size of a generating set of $\SO_5(\BZ[{1}/{5}])$ then
any $n$ elements in $\SO_5(\BR)$ can be slightly deformed to a
generating set of the group $\SO_5(\BZ[{1}/{5}])$ (recall that
$\SO_5(\BZ[{1}/{5}])$ is dense in $\SO_5(\BR)$ and has Kazhdan property (T), c.f. \cite{margulis}). Similarly, any
$5$ elements in $\SO_3(\BR)$ can be deformed to a generating set
of a surface group.

\medskip
Let $G$ be a connected compact Lie group. Since the topology of
$G$ is metrizable, we may assume that it is endowed with a metric
$d$, and by averaging over all left and right translations with
respect to the Haar measure, we may assume that $d$ is left and
right invariant. Let $\gC\leq G$ be an $(n-1)$--generated dense
subgroup of $G$, let $s_1,\ldots,s_n$ be arbitrary $n$ elements of
$G$, and let $\gep>0$. We will explain how to perturb the $s_i$ to
some $t_i$ which satisfy $d(t_i,s_i)<\gep,~i=1,\ldots ,n$ and
$\langle t_1,\ldots,t_n\rangle=\gC$.

We will first treat the case where $G$ is semisimple, since our
proof in this case gives more than what stated in \ref{thm1}.

\subsection{The semisimple case}
Let $G$ be a connected compact semisimple Lie group. It is well
known that $G$ admits a structure of an algebraic group over
$\BR$, and hence is equipped with a Zariski topology. Let
$\mathfrak{g}$ denote the Lie algebra of $G$, write
$\mathfrak{g}=\bigoplus_{i=1}^p\mathfrak{g}_i$ where
$\mathfrak{g}_i$ are the simple factors of $\mathfrak{g}$, and set
$\mathcal{A}=\bigoplus_{i=1}^p\text{End}(\mathfrak{g}_i)$. Let
$G_i$ be the subgroup of $G$ corresponding to $\mathfrak{g}_i$;
recall that $G$ is an almost direct product of the $G_i$. Let
$\overline{G}_i$ be the quotient of $G$ by the product of the
$G_j,~j\neq i$ and let $\pi_i:G\to\overline{G}_i$ be the canonical
projection. For each $i$, the restriction of the Adjoint
representation to $G_i$ is irreducible on $\mathfrak{g}_i$ and
hence, by Burnside's lemma,
$\text{span}(\text{Ad}(G))=\mathcal{A}$; It is well known that
$\mathcal{A}$ is generated by two elements which can be taken from
$\text{Ad}(G)$, for instance one can take two elements in
$\text{Ad}(G)$ which generate a dense subgroup (c.f.
\cite{kuranishi}).

\begin{lem}
Let $m\in\BN$ and let $g_1,\ldots,g_m$ be arbitrary $m$ elements
in $G$. The set
$$
 \gO (g_1,\ldots,g_m):=\{ g\in G:
 \text{Ad}(g),\text{Ad}(g_1),\dots,
 \text{Ad}(g_m)~\text{~generate~}~\mathcal{A}\}
$$
is Zariski open in $G$.
\end{lem}

\begin{proof}
If $\gO (g_1,\ldots,g_m)$ is empty there is nothing to prove.
Assume that $\gO (g_1,\ldots,g_m)\neq\emptyset$ and let $g_0\in\gO
(g_1,\ldots,g_m)$. Let $d=\dim(\mathcal{A})$. It follows from the
definition of $\gO (g_1,\ldots,g_m)$ that there are $d$ words in
$m+1$ letters $W_1,\ldots,W_{d}$ which, when evaluated at the
point $(\text{Ad}(g_0),\text{Ad}(g_1),\ldots,\text{Ad}(g_m))$,
span $\mathcal{A}$. The set
$$
 \{ g\in G:\text{span}\{ W_i\big(\text{Ad}(g),\text{Ad}(g_1),\ldots,\text{Ad}(g_m)\big): i=1,\ldots,d
 \}=\mathcal{A}\}
$$
is then a Zariski open subset of $\gO (g_1,\ldots,g_m)$ which
contains $g_0$. Since $g_0$ is an arbitrary element of $\gO
(g_1,\ldots,g_m)$ it follows that $\gO (g_1,\ldots,g_m)$ is
Zariski open.
\end{proof}

Suppose now that $g_1,\ldots,g_k\in G$ are such that
$\text{Ad}(g_1),\ldots,\text{Ad}(g_k)$ generate $\mathcal{A}$,
then for each $1\leq i\leq k$, the set\footnote{The hat above the
$i$'th element means that we exclude it.}
$\gO(g_1,\ldots,\hat{g_i},\ldots,g_{k})$ is Zariski open, and
since it contains $g_i$ it is nonempty and hence Zariski dense in
$G$ as $G$ is Zariski connected. It follows that also the set
$$
 \ti\gO(g_1,\ldots,g_k):=\cap_{i=1}^k\gO(g_1,\ldots,\hat{g_i}\ldots,g_{k})
$$
is Zariski open and nonempty, hence dense in $G$ with respect to
the housdorf topology. Moreover, we have:

\begin{lem}\label{lem:density}
Let $g_1,\ldots,g_k\in G$ be $k$ elements such that
$\text{Ad}(g_1),\ldots,\text{Ad}(g_k)$ generate $\mathcal{A}$. The
projection of $\langle g_1,\ldots,g_k\rangle$ to each simple
factor $\overline{G_i},~i=1,\ldots,p$ of $G$ is either finite or
dense. If for every $i\leq p$ this projection is infinite, then
$\langle g_1,\ldots,g_k\rangle$ is dense in $G$.
\end{lem}

\begin{proof}
Let $C$ be the closure of $\langle g_1,\ldots,g_k\rangle$ and
$C^\circ$ its identity connected component. Then $C^\circ$ is a
connected compact Lie subgroup of $G$ and $|C/C^\circ|<\infty$.
Moreover, since the Lie algebra of $C^\circ$ is stable under the
elements $\text{Ad}(g_i),i=1,\ldots,k$ which generate
$\mathcal{A}$, $C^\circ$ is normal in $G$. Since for
$i=1,\ldots,p$, $\overline{G_i}$ is almost simple, it follows that
the compact connected normal subgroup $\pi_i(C^\circ)$ is either
trivial or equal to $\overline{G}_i$, which means that
$\pi_i(\langle g_1,\ldots,g_k\rangle)$ is either finite or dense.
Finally, if $\pi_i(\langle g_1,\ldots,g_k\rangle)$ is infinite for
each $i\leq p$ then the normal subgroup $C^\circ$ projects onto
every simple factor of $G$, and hence $C^\circ=G$, i.e. $\langle
g_1,\ldots,g_k\rangle$ is dense.
\end{proof}

We will also use the following well known:

\begin{lem}\label{lem:U-open-dense}
The set $U$ of all pairs $(x_1,x_2)\in G\times G$ for which the
group $\langle x_1,x_2\rangle$ is dense in $G$, is open, dense and
of full Haar measure in $G\times G$.
\end{lem}

\begin{proof}
The second and third properties follow from the fact that $U$
contains the the following set which is clearly dense and has full
Haar measure:
$$
 \{ (x_1,x_2)\in G\times G:\text{Ad}(x_1),\text{Ad}(x_2)~
 \text{~generate~}~\mathcal{A},~\text{~and~}~\pi_j(x_1)~\text{~is non-torsion~}~\forall j\leq p
 \}.
$$
It is well known that $U$ contains an open subset $V$ near the
identity of $G\times G$ (c.f. \cite{GZ} or \cite{BG}). To see that
$U$ is open, note that if $(x_1,x_2)\in U$ then, as $\langle
x_1,x_2\rangle$ is dense, there are two words in two letters
$W_1,W_2$ such that $\big( W_1(x_1,x_2),W_2(x_1,x_2)\big)$ belongs
to $V$. It follows that if $U_i,~i=1,2$ are sufficiently small
neighborhoods of $x_i$ in $G$ then
($(W_1(y_1,y_2),W_2(y_1,y_2))\in V$ for any $(y_1,y_2)\in
U_1\times U_2$, which implies that $\langle y_1,y_2\rangle$ is
dense in $G$. Therefore $U_1\times U_2\subset U$, and since
$(x_1,x_2)$ is arbitrary, $U$ is open.
\end{proof}

In order to prove Theorem \ref{thm1} we will use repeatedly the so called product replacement moves  
which allows replacing one generating set
$(\gc_1,\ldots,\gc_k)$ for the group $\gC$ by another generating
set of the same cardinality by multiplying one $\gc_i$ by some
$\gc_j^{\pm 1}$ where $j\neq i$ or more generally by an element of
the group $\langle \gc_j:j\neq i\rangle$. These operations are
also called Nielsen transformations.

Let $(\gc_1,\ldots,\gc_{n-1})$ be a generating set for $\gC$.
Applying Selberg's lemma to the projection of $\gC$ to any simple
factor $\overline{G}_i$ of $G$ and taking the intersection, we see
that $\gC$ contains a subgroup of finite index $\gC_0$ such that
the projection of any element of $\gC_0$ to any simple factor of
$G$ is either trivial or non-torsion. Since $\gC$ is dense and $G$
is connected, also $\gC_0$ is dense. Pick
$\gc_n\in\gC_0\cap\ti\gO(\gc_1,\ldots,\gc_{n-1})$ with
$\pi_j(\gc_n)\neq 1,\forall j\leq p$. Then from Lemma
\ref{lem:density} and the assumption that $\gC$ is dense we get
that $\langle\gc_1,\ldots,\hat\gc_i,\ldots,\gc_n\rangle$ is dense
in $G$ for every $1\leq i\leq n$. Fix two open sets $U_i\subset
B_\gep (s_i),~i=1,2$ such that $U_1\times U_2\subset U$, and
$U_1\subset\ti\gO (\gc_2,\ldots,\gc_n)$, where $U\subset G\times
G$ is the set defined in Lemma \ref{lem:U-open-dense}. Since
$\langle\gc_2,\ldots,\gc_n\rangle$ is dense in $G$ there is some
$\hat\gc_1\in\langle\gc_2,\ldots,\gc_n\rangle$ for which
$t_1:=\hat\gc_1\gc_1\in U_1$. Since $t_1$ belongs to $\ti\gO
(\gc_2,\ldots,\gc_n)$, it follows from Lemma \ref{lem:density} and
the property of $\gc_n$ that the group $\langle
t_1,\gc_3,\ldots,\gc_n\rangle$ is dense in $G$. Chose
$\hat\gc_2\in\langle t_1,\gc_3,\ldots,\gc_n\rangle$ such that
$t_2=\hat\gc_2\gc_2$ lies in $U_2$. It follows that $\langle
t_1,t_2\rangle$ is dense in $G$ and we can pick
$\hat\gc_i\in\langle t_1,t_2\rangle,~i=3,\ldots,n$ so that
$t_i:=\hat\gc_i\gc_i\in B_\gep (s_i)$. This completes the proof of
Theorem \ref{thm1} in the semisimple case.\qed

\begin{rems}
(1) The argument above, with some simple modifications, can be
applied also to dense subgroups $\gC$ of the form
$\gC=\langle\gD,\gc\rangle$ where $\gD$ is an $(n-1)$--generated
dense subgroup of $G$, and $\gc\in G\setminus\gD$. It is therefore
natural to ask wether any $n$--generated dense subgroup of $G$ is
of that form. This question makes sense for every $n\geq 3$, and the answer may depend on $n$ and $G$.
Nir Avni pointed out to me that for every $G$ there is some integer $n(G)$ such that the answer is affirmative for all $n\geq n(G)$. It is of interest wether $n(G)$ is always $3$ or some other constant independent of $G$. 
For $n\geq n(G)$ one can assume that $\gC$ is $n$--generated in Theorem \ref{thm1}.

(2) It was shown in \cite{GZ} that any $k$--generated dense
subgroup $\gC$ of a connected semisimple Lie group $G$ admits a
generating set of cardinality $(k+2)$ which lies arbitrarily close
to the identity of $G$, and hence there is no uniform Kazhdan
constant for generating sets of size $k+2$ for the natural
representation of $\gC$ on $L_2(G)$. It follows from Theorem \ref{thm1} that
when $G$ is compact, one can replace $k+2$ by $k+1$. In view of the previous remark, for $k\geq n(G)$ one can even omit the "$+1$".

(3) The argument above can be also applied for some non connected
compact groups. M. Abert and L. Pyber have shown to me that using
results from \cite{Dunwoody} one can prove the following: Let $G$
be a topologically $k$--generated pro-finite-soluble group and
$\Gamma$ a $(k+1)$--generated dense subgroup of $G$. Then for any
$k+1$ elements $s_1,\ldots,s_{k+1}$ which topologically generate
$G$ and any open normal subgroup $N\lhd G$, one can choose $t_i\in
s_iN,~i=1,\ldots,k+1$ such that $\Gamma=\langle
t_1,\ldots,t_{k+1}\rangle$. Moreover any $k$--generated profinite
group has a $k+1$--generated dense subgroup which is not virtually
free.

(4) For non-compact Lie groups E. Ghys asked the following related
question: Let $H$ be a connected non-compact simple Lie group.
Does every $n$--tuple $(s_1,\ldots,s_n)$ which generates a dense
subgroup of $H$ admit an arbitrarily small deformation
$(t_1,\ldots,t_n)$ for which the group $\langle
t_1,\ldots,t_n\rangle$ is not free? We refer to \cite{Soifer} for
other related problems.
\end{rems}


\subsection{The ergodicity of $\text{Aut}(F_n)$ on $\text{Hom}(F_n,G)$}

By fixing a free generating set for the free group $F_n$ we
may identify the deformation variety $\text{Hom}(F_n,G)$ with
$G^n$. The automorphism group $\text{Aut}(F_n)$ acts on
$\text{Hom}(F_n,G)$ by pre-compositions. The Nielsen operations on
generating sets (the product replacement moves) when applied to
the generators of the free group $F_n$ correspond to elements of
$\text{Aut}(F_n)$. Let us denote by $L_{i,j}$ the operation of
replacing the $j$'th generator by its product from the left with
the $i$'th generator, e.g.
$L_{2,1}(g_1,g_2,\ldots,g_n)=(g_2g_1,g_2,\ldots,g_n)$. Arguing as
above one can show that the $\text{Aut}(F_n)$--orbit of almost any
point in $G^n$ is dense. More precisely, let $Y\subset G^n$ be the
set of $n$--tuples consisting of $n$ elements, such that the
Adjoint of any $n-1$ of them generate the algebra $\mathcal{A}$
and the projection of each to any simple factor of $G$ is
non-torsion, and let $Y_0=\cap_{\gs\in\text{Aut}(F_n)}\gs (Y)$.
Then $Y_0$ is $\text{Aut}(F_n)$--invariant, has full measure and
the orbit of any element in $Y_0$ is dense in $G^n$. Furthermore,
we have:

\begin{thm}\label{thm:ergodic}
Let $G$ be a compact connected semisimple Lie group, and let
$n\geq 3$. The action of $\text{Aut}(F_n)$ on $G^n$ is ergodic.
\end{thm}

\begin{proof}
Assume the contrary, and let $A\subset G^n=G\times
G\times\ldots\times G$ be an $\text{Aut}(F_n)$ almost invariant
measurable subset which is neither null nor conull. Since
$\text{Aut}(F_n)$ is countable we may assume that $A$ is invariant
rather than almost invariant. Now since the action of $G$ on
itself by left translations is clearly ergodic, it follows from
our assumption that $A$ is not (almost) invariant under the action
by left translations of one of the $n$ factors, say the first, of
$G^n$. Thus by Fubini's theorem, for a set of positive measure of
$(g_2,\ldots,g_n)\in G^{n-1}$ we have that $\{ g\in
G:(g,g_2,\ldots,g_n)\in A\}$ is neither null nor conull in $G$.
Let us fix a point $(g_2,\ldots,g_n)$ in this subset such that the
pair consisting of the first two components $(g_2,g_3)$ belongs to
the set of full measure $U\subset G\times G$ introduced in Lemma
\ref{lem:U-open-dense}. Note that the orbits of the action of
$\langle g_2,g_3\rangle$ by left translations on $G$ coincides
with (the projection to the first factor of) the orbits of
$\langle L_{2,1},L_{3,1}\rangle$ on $\{(g,g_2,\ldots,g_n):g\in
G\}$. Set
$$
 A_1:=\{ g\in G:(g,g_2,\ldots,g_n)\in A\}.
$$
By our assumption, $A_1$ is neither null or conull. This however
is a contradiction since the group $\langle g_2,g_3\rangle$ is
dense in $G$ and hence acts ergodically on $G$.
\end{proof}

Theorem \ref{thm:ergodic} clearly implies, and is actually equivalent to
(see \cite{Goldman} Lemma 3.1):

\begin{thm}\label{thm:erg-2}
Let $G$ be a compact connected semisimple Lie group and let $n\geq
3$. The action of $\text{Out}(F_n)$ on $\text{Hom}(F_n,G)/G$ is
ergodic.
\end{thm}

Theorem \ref{thm:erg-2} was conjectured by W.M. Goldman in
\cite{Goldman}, and was proved there under the assumption that
$G=\text{SU}(2)$. As explained in the last paragraph of
\cite{Goldman}, the general case follows from the semisimple one,
i.e. one can omit the assumption that the connected compact Lie
group $G$ is semisimple in Theorems \ref{thm:ergodic} and
\ref{thm:erg-2}. As pointed out in \cite{Goldman} the
assumption that $n\geq 3$ is necessary, since the function $\text{trace}([x,y])$ is $\text{Aut}(F_2)$ invariant on $G^2$.

Applying the conclusions of Theorems \ref{thm:ergodic} and \ref{thm:erg-2} to $G^2$ instead of $G$ we obtain furthermore:

\begin{thm}
Let $G$ be a compact connected Lie group and let $n\ge 3$. Then the action of $\text{Aut}(F_n)$ on $\text{Hom}(F_n,G)$, and the action of $\text{Out}(F_n)$ of $\text{Hom}(F_n,G)/G$ are weakly mixing.
\end{thm}

\begin{proof}
Since $\text{Hom}(F_n,G)\times\text{Hom}(F_n,G)$ is canonically isomorphic to $\text{Hom}(F_n,G\times G)$, and since $\text{Hom}(F_n,G)/G\times\text{Hom}(F_n,G)/G$ is canonically isomorphic to $\text{Hom}(F_n,G\times G)/G\times G$, the assertions follows by applying Theorems \ref{thm:ergodic} and \ref{thm:erg-2} to the connected compact Lie group $G\times G$.
\end{proof}

\begin{rem}
An $(n\ge 2)$--tuple $(g_1,\ldots,g_n)\in G^n$ is said to have a spectral gap if the maximal eigenvalue of the corresponding Hecke operator on $L_2(G)$ is isolated in the spectrum, or equivalently, if the group $\langle g_1,\ldots,g_n\rangle$ acts strongly ergodically, by left multiplications, on $G$. D. Fisher \cite{Fisher} showed that Goldman's theorem about the ergodicity of the $\text{Aut}(F_n)$ action on $\SO(3)^n$ implies that for $n\ge 3$ the set of $n$--tuples in $\SO(3)$ which posses a spectral gap is either null or conull. Indeed, this set is measurable and $\text{Aut}(F_n)$--invariant. Theorem \ref{thm:ergodic} implies that the same conclusion holds when $\SO(3)$ is replaced by any compact connected Lie group $G$. Furthermore, Theorem \ref{thm1} implies that if $G$ admits one $(n\ge 2)$--tuple with a spectral gap, then the set of $(n+1)$--tuples with spectral gap is dense in $G^{n+1}$.
\end{rem}


\subsection{The proof of Theorem \ref{thm1} in the case where $G$ is not necessarily semisimple}
Let us come back to the proof of Theorem \ref{thm1} in the general
case. Let $G$ be an arbitrary connected compact Lie group. Then
$G$ is an almost direct product of $G'$ and $Z$ where $G'$ is the
commutator group of $G$ and is semisimple, and $Z$ is the center
of $G$. We do not assume that $G'$ and $Z$ are non-trivial. Let
$D\subset G\times G$ be the set of all pairs which generate a
dense subgroup of $G$.

\begin{lem}
$D$ is dense in $G\times G$.
\end{lem}

\begin{proof}
Let $V\subset G\times G$ be the pre-image of the open dense subset
of $G/Z\times G/Z$ that one gets from Lemma \ref{lem:U-open-dense}
applied to the semisimple group $G/Z$ (if $G=Z$ take $V=G\times
G$). If $(x_1,x_2)\in V$ then $\langle x_1,x_2\rangle\cap G'$ is
dense in $G'$.

The quotient $G/G'$, being a connected compact abelian Lie group,
is isomorphic to a product of circle groups, and it is well known
that in such a group the set of elements which generate a dense
cyclic subgroup is dense. In fact, an element in $U(1)^m$
generates a dense cyclic subgroup iff its normalized coordinates
are independent together with $1$ over the rational. Let $\pi
:G\to G/G'$ denote the canonical projection. The set
$$
 \{ (x_1,x_2)\in V:\langle\pi (x_1)\rangle~\text{~is
 dense in~}~G/G'\}
$$
is clearly dense in $G\times G$ and contained in $D$.
\end{proof}

Fix a pair $(g_1,g_n)\in D\cap (B_\frac{\gep}{2}(s_1)\times
B_\frac{\gep}{2}(s_n))$, and choose finitely many words
$W_j,~j=1,\ldots,k$ in two letters, such that the set $\{
W_j(g_1,g_n):j=1,\ldots,k\}$ form an $\gep/2$--net in $G$. Then
fix $\gep_1<\gep/2$ sufficiently small so that for any
$(x_1,x_n)\in B_{\gep_1}(g_1)\times B_{\gep_1}(g_n)$ the set $\{
W_j(x_1,x_n):j=1,\ldots,k\}$ is an $\gep$--net in $G$. Let
$\{\gc_1,\ldots,\gc_{n-1}\}$ be a generating set for the dense
subgroup $\gC$. We will use the following:

\begin{clm}
The set $X=\{ x\in G: \langle
\gc_2,\ldots,\gc_{n-1},x\rangle~\text{~is dense in~}~ G\}$ is
dense in $G$.
\end{clm}

Indeed, $X$ contains every $x$ which satisfies the following
conditions:
\begin{enumerate}
\item
$\text{Ad}(\gc_2),\ldots,\text{Ad}(\gc_{n-1}),\text{Ad}(x)~\text{~generate
the algebra span(Ad}(G))$.

\item The projection of $x$ to every simple factor of $G$ is
non-torsion.

\item The projection of $x$ to $G/G'$ generates a dense subgroup
of $G/G'$.
\end{enumerate}
Condition (1) defines an open dense subset of $G$, while the set
of elements which satisfies (2) and (3) is clearly dense.

Fix $x_n\in X\cap B_{\gep_1}(g_n)$ and let $W$ be a word in $n-1$
letters such that $W(\gc_2,\ldots,\gc_{n-1},x_n)\in
B_{\gep_1}(g_1)\gc_1^{-1}$. Since $\gC$ is dense in $G$ we can
pick $t_n\in\gC\cap B_{\gep_1}(g_n)$ sufficiently close to $x_n$
so that $W(\gc_2,\ldots,\gc_{n-1},t_n)$ still lies in
$B_{\gep_1}(g_1)\gc_1^{-1}$. Then
$t_1:=W(\gc_2,\ldots,\gc_{n-1},t_n)\gc_1$ lies in
$B_{\gep_1}(g_1)$. It follows that the set
$$
 \{ W_j(t_1,t_n):j=1,\ldots, k\}
$$
is an $\gep$--net in $G$. Since the metric on $G$ is invariant
under right translations, there are $n-2$ (not necessarily
distinct) elements of this $\gep$--net:
$W_{j_i}(t_1,t_n),i=2,\ldots,n-1$ such that the elements
$t_i:=W_{j_i}(t_1,t_n)\gc_i,i=2,\ldots,n-1$ belong to $B_\gep
(s_i)$ respectively. Since $(t_1,\ldots,t_n)$ was obtained from
$(\gc_1,\ldots,\gc_{n-1},t_n)$ by product replacement moves,
$(t_1,\ldots,t_n)$ generates $\gC$ and the theorem is proved.\qed


\subsection{A concrete example}\label{example1}
Let $G$ be a connected compact Lie group with
center $Z$, and let $a,b,c\in G$ be arbitrary three elements.
Slightly deforming $a$ and $b$ to $a',b'$ we may assume that:
\begin{enumerate}
\item $a'$ is regular and non-torsion.

\item The group $\langle a'Z,b'Z\rangle$ is dense in $G/Z$.
\end{enumerate}

The centralizer $A=Z_G(a')$ of $a'$ is a maximal torus in $G$.
Since all maximal tori of $G$ are conjugate and their union is
equal to $G$, we can find $\gc\in\langle a',b'\rangle$ such that
$\gc A\gc^{-1}$ passes arbitrarily close to $c$. Then we can
slightly deform $c$ to some $c'\in\gc A\gc^{-1}$ such that $c'$
and $\gc a'\gc^{-1}$ are independent over $\BZ$, i.e. $\langle c',
\gc a'\gc^{-1}\rangle\cong\BZ^2$. The group $\langle
a',b',c'\rangle$ is not virtually free since it contains a copy of
$\BZ^2$.


\section{The case $n=2$}

The argument above, using the product replacement method, does not
apply for $n=2$. In this case we prove the following:

\begin{thm}\label{thm2}
Let $G$ be a compact connected semisimple Lie group, and let
$a,b\in G$. There is an arbitrarily small deformation $a',b'$ of
$a,b$ such that $\langle a',b'\rangle$ is dense\footnote{By
replacing the word ``dense" with the word ``infinite" in the
conclusion of the theorem, one may assume that $G$ is only
non--abelian rather than semisimple.} in $G$ and has Serre's
property (FA). In particular $\langle a',b'\rangle$ is not
virtually free.
\end{thm}

Recall that a group $\gC$ has Serre's property (FA) if every action
of $\gC$ on a tree admits a global fixed point.

A finitely generated infinite virtually non-abelian-free group,
being quasi isometric to $F_2$, has infinitely many ends, and
hence by Stalling's theorem splits over a finite group, and
therefore acts minimally on some tree with finite edge
stabilizers. Additionally, it is well known that an infinite
virtually cyclic group admits a transitive action on the linear
tree. Therefore, an infinite group with property (FA) is not
virtually free.

Yves de Cornulier suggested to use the following result of Serre
\cite{Serre}:

\begin{lem}[Serre \cite{Serre}]\label{lem:serre} Let $\gC=\langle a,b\rangle$
and suppose that $a,b$ and $ab$ are torsion. Then $\gC$ has property (FA).
\end{lem}

For the sake of completeness let us give a proof for Lemma
\ref{lem:serre}: Suppose that $\gC$ acts on a tree $T$, and let
$T_a,T_b$ be the fixed subtrees of $a,b$, which are non-empty
since $a$ and $b$ are torsion. Then $T_a\cap T_b$ is non-empty,
for otherwise the segment $I$ connecting $T_a$ to $T_b$ would be
oriented with $(ab)\cdot I$ which in turn would imply that $ab$ is
hyperbolic, in contrary to the assumption that it is torsion.

\medskip

We will show that any two elements can be deformed to elements
which satisfy the conditions of Lemma \ref{lem:serre} and generate
a dense subgroup. For the sake of simplicity let us assume that
$G$ is simple. The semisimple case follows easily by considering
simultaneously all simple factors of $G$. By Lemma
\ref{lem:U-open-dense} the set of couples which generate a dense
subgroup in $G$ is open and dense in $G\times G$. Start with
arbitrary $a,b\in G$. Slightly deforming them we may assume that

\begin{enumerate}
\item Both $a$ and $b$ are torsion.

\item Both $a$ and $b$ are regular, in the sense that their
centralizers are maximal tori in $G$.

\item The group $\langle a,b\rangle$ is dense in $G$.
\end{enumerate}

We will show that it is possible to deform $a$ and $b$ within
their conjugacy class, i.e. find $g,h\in G$ close to the identity,
such that $a^gb^h=(gag^{-1})(hbh^{-1})$ would be a torsion
element.

Denote by $\mathcal{T}_a=Z_G(a)$ the centralizing torus of $a$, by $\mathcal{T}_b$ the
one of $b$ and by $\mathfrak{t}_a,\mathfrak{t}_b$ the
corresponding abelian Lie subalgebras of
$\mathfrak{g}=\text{Lie}(G)$.

\begin{clm}\label{clm}
$\mathfrak{t}_a\cap\mathfrak{t}_b=\{ 0\}$.
\end{clm}

Since $\mathfrak{t}_a\cap\mathfrak{t}_b$ is in the null space of
both $\text{Ad}(a)$ and $\text{Ad}(b)$ while $\langle
a,b\rangle$ is dense, the claim follows as the center of
$\mathfrak{g}$ is trivial.

\medskip

For $x\in G$, let $C(x)=\{ gxg^{-1}:g\in G\}$ denote the conjugacy
class of $x$. We will derive Theorem \ref{thm2} from the following
proposition, in which we do not require that $G$ is compact:

\begin{prop}\label{prop:pco}
Let $G$ be a connected simple Lie group and $a,b\in G$ two regular
elements in general position, i.e.
$\mathfrak{t}_a\cap\mathfrak{t}_b=\{ 0\}$. The restriction of the
product map $G\times G\to G$ to the conjugacy class $C(a)\times
C(b)$ of $(a,b)$ is open in a neighborhood of $(a,b)$.
\end{prop}

\begin{proof}
First note that $\mathfrak{t}_a=\ker (\text{Ad}(a)-1)$ and
$\mathfrak{t}_b=\ker (\text{Ad}(b)-1)$. Now
$$
 (\text{Ad}(x^{-1})-1)(\mathfrak{g})=\big(\ker
 (\text{Ad}(x)-1)\big)^\bot,
$$
where the orthogonal complement is taken with respect to the
Killing form, which is non-degenerate as $G$ is simple. Since
$\mathfrak{t}_a\cap\mathfrak{t}_b=\{ 0\}$ it follows that
\begin{equation}\label{eq1}
 (\text{Ad}(a^{-1})-1)(\mathfrak{g})+(\text{Ad}(b^{-1})-1)(\mathfrak{g})=\mathfrak{g}.
\end{equation}

Since $G$ admits a faithful linear representation, we my assume
that it is linear. Let us compute the tangent space to the
conjugacy class $C(a)$ at $a$, identified via the Lie algebra of
left invariant vector fields as a subspace of the tangent space
$\mathfrak{g}=T_1(G)$. For $X\in\mathfrak{g}$ we have
$$
 \frac{d}{dt}\big|_{t=0}(\exp (tX)a\exp (-tX))=Xa-aX,
$$
and by multiplying from the left by $a^{-1}$ we obtain:
$$
 L_{a^{-1}}T_a(C(a))=(\text{Ad}(a^{-1})-1)(\mathfrak{g}),
$$
where $T_a(C(a))$ denotes the tangent space of $C(a)$ at $a$ viewed as a subspace of $T_a(G)$.
Similarly,
$L_{b^{-1}}T_b(C(b))=(\text{Ad}(b^{-1})-1)(\mathfrak{g})$.

Finally, the differential of the product map $G\times G\to G$ at
$(a,b)$, again identified via a left translation as a subspace of
$T_{(1,1)}(G\times G)=\mathfrak{g}\oplus\mathfrak{g}$, evaluated
at $(X,Y)$ is easily seen to be $\text{Ad}(b^{-1})(X)+Y$. It
follows that the image of the differential at $(a,b)$ of the
product map $C(a)\times C(b)\to G$ is:
$$
 \text{Ad}(b^{-1})(\text{Ad}(a^{-1})-1)(\mathfrak{g})+(\text{Ad}(b^{-1})-1)(\mathfrak{g}).
$$
Since the second summand is $\text{Ad}(b)$--invariant, we derive
from (\ref{eq1}) that this differential is onto. Therefore, the
proposition follows from the implicit function theorem.
\end{proof}

Now since the torsion elements are dense in $G$, we can pick, by
Proposition \ref{prop:pco} and Lemma \ref{lem:U-open-dense},
$g,h\in G$ arbitrarily close to $1$ such that $\langle
a^g,b^h\rangle$ is still dense and $a^gb^h$ is torsion. This
completes the proof of Theorem \ref{thm2}.\qed

\medskip

\noindent {\bf Acknowledgement:} I would like to thank Gregory
Soifer, Jean Francois Quint, Emmanuel Breuillard, Miklos Abert,
Laszlo Pyber, Nir Avni, Shahar Mozes, Alireza Salehi Golsefidy and Uri Bader
for helpful discussions. I also thank Yves de Cornulier for his helpful suggestion to use Lemma \ref{lem:serre} of Serre.
I am very grateful to Hee Oh for inviting me to Caltech in April 2006 where I had the main ideas behind this work.

\end{document}